\documentclass[12pt,reqno]{amsart}
\usepackage{amsmath}
\usepackage{amssymb}
\usepackage{amsthm}
\usepackage{latexsym}
\usepackage{mathrsfs}

\newtheorem{theorem}{Theorem}[section]
\newtheorem{lemma}{Lemma}[section]
\newtheorem{corollary}{Corollary}[section]
\newtheorem{remark}{Remark}[section]
\newtheorem{definition}{Definition}[section]
\newtheorem{proposition}{Proposition}[section]
\newtheorem{example}{Example}[section]
\newtheorem{assumption}{Assumption}[section]
\numberwithin{equation}{section}
\newcommand{\bth}{\begin{theorem}}
\newcommand{\ethe}{\end{theorem}}

\newcommand{\bre}{\begin{remark}}
\newcommand{\ere}{\end{remark}}

\newcommand{\ble}{\begin{lemma}}
\newcommand{\ele}{\end{lemma}}
\newcommand{\bde}{\begin{definition}}
\newcommand{\ede}{\end{definition}}

\newcommand{\bco}{\begin{corollary}}
\newcommand{\eco}{\end{corollary}}

\newcommand{\bpr}{\begin{proposition}}
\newcommand{\epr}{\end{proposition}}

\newcommand{\bexer}{\begin{exercise}}
\newcommand{\eexer}{\end{exercise}}
\newcommand{\breh}{\begin{hint}}
\newcommand{\ereh}{\end{hint}}

\newcommand{\halmos}{\hfill \qed}

\newcommand{\bexam}{\begin{example}}
\newcommand{\eexam}{\end{example}}

\newcommand{\pr} {{\bf Proof.}}

\newcommand{\bfi}{\begin{fig}}
\newcommand{\efi}{\end{fig}}

\newcommand{\beao}{\begin{eqnarray*}}
\newcommand{\eeao}{\end{eqnarray*}\noindent}

\newcommand{\beam}{\begin{eqnarray}}
\newcommand{\eeam}{\end{eqnarray}\noindent}
\newcommand{\E}{\mathbf{E}}
\newcommand{\PP}{\mathbf{P}}

\newcommand{\xto}{x\to\infty}

\newcommand{\bF}{\overline{F}}

\newcommand{\bV}{\overline{V}}
\newcommand{\bA}{\overline{A}}

\newcommand{\bbr}{{\mathbb R}}

\newcommand{\bbb}{{\mathbb B}}

\newcommand{\bbn}{{\mathbb N}}

\newcommand{\vep}{\varepsilon}

\usepackage{soul}
\setstcolor{red}

\usepackage{color}

\begin{document}
\title[Uniform asymptotic estimates for ruin probabilities]{Uniform asymptotic estimates for ruin probabilities of a multidimensional risk model with  c\`{a}dl\`{a}g returns and multivariate heavy tailed claims}

\author[D.G. Konstantinides, C.D. Passalidis]{Dimitrios G. Konstantinides,\,Charalampos  D. Passalidis} 

\address{Dept. of Statistics and Actuarial-Financial Mathematics,
University of the Aegean,
Karlovassi, GR-83 200 Samos, Greece}
\email{konstant@aegean.gr,\;sasm23002@sas.aegean.gr.}

\date{{\small \today}}

\begin{abstract} 
We study a multidimensional renewal risk model, with common counting process and  c\`{a}dl\`{a}g returns. Considering that the claim vectors have common distribution from some multivariate distribution class with heavy tail, are mutually weakly dependent, and each one has arbitrarily dependent components, we obtain uniformly asymptotic estimations for the probability of entrance of discounted aggregate claims into a some rare sets, over a finite time horizon. Direct consequence of the claim behavior is the estimation of the ruin probability of the model in some ruin sets. Further, restricting the distribution class of the claim vectors in the multivariate regular variation, the estimations still hold uniformly over the whole time horizon.   
\end{abstract}

\maketitle
\vspace{3mm}


\section{Introduction} \label{sec.KMP.1}

We consider an insurer, who operates $d$-business lines simultaneously, with $d \in \bbn$, that follow a common, renewal, claim counting process. Because of the intense competition, the modern insurance companies invest the surplus either in risk free or in risky assets. In our model, we assume that the price process of the investment portfolio, is described by the stochastic process $\{e^{\xi(t)}\,,\;t\geq 0\}$, where $\{\xi(t)\,,\;t\geq 0\}$, represents a c\`{a}dl\`{a}g stochastic process, with $\xi(0)=0$, which is general enough, found in the literature of uni-variate or multivariate risk models, since it contains the constant interest rate, the zero interest rate or any L\'{e}vy process, see for example \cite {guo:wang:2013}, \cite{paulsen:2002}, \cite{paulsen:2008}, \cite{heyde:wang:2009}, \cite{kalashnikov:norberg:2002}, etc. Furthermore, we accept that the insurer receives premiums, whose density is depicted by the vector ${\bf c}(t)=(c_1(t),\,\ldots,\,c_d(t))$, with $t \geq 0$, and $c_i(t)$ is the density function of premiums in $i$-th line of business, for $i=1,\,\ldots,\,d$. In what follows, the density of premiums is bounded, namely $0\leq c_i(t) \leq M_i$, for some constants $M_i\geq 0$ and for any $t \geq 0$ and $i=1,\,\ldots,\,d$.

Also we consider the total initial capital of the insurer is equal to $x$ and it is allocated in the $d$ lines of business according to weights $l_1,\,\ldots,\,l_d > 0$, such that $l_1+\cdots + l_d =1$. Finally, the $i$-th claim vector ${\bf X}^{(i)}=(X_1^{(i)},\,\ldots,\,X_d^{(i)})$, can allow zero components, but not all of them, and it arrives at time moment $T_i$, with $i \in \bbn$ and conventionally we put $T_0=0$. The $\{T_i\,,\;i\in \bbn\}$ constitute a renewal counting process, symbolically $\{ N(t)\,,\;t\geq 0\}$, with finite renewal function $\lambda(t) =\E[N(t)]=\sum_{i=1}^\infty \PP[T_i\leq t]$.

The discount surplus process of the insurer, can be written through the relations
\beam \label{eq.KMP.1.1}
\left( 
\begin{array}{c}
U_{1}(t) \\ 
\vdots \\ 
U_{d}(t) 
\end{array} 
\right) =x\,\left( 
\begin{array}{c}
l_{1} \\ 
\vdots \\ 
l_{d} 
\end{array} 
\right) +\left( 
\begin{array}{c}
\int_{0-}^{t}e^{-\xi(s)}\,c_1(s)\,ds \\ 
\vdots \\ 
\int_{0-}^{t}e^{-\xi(s)}\,c_d(s)\,ds 
\end{array} 
\right) -\left( 
\begin{array}{c}
\sum_{i=1}^{N(t)} X_{1}^{(i)}\,e^{-\xi(T_{i}) } \\ 
\vdots \\ 
\sum_{i=1}^{N(t)} X_{d}^{(i)}\,e^{-\xi(T_{i}) } 
\end{array} 
\right)\,,
\eeam 
for any $t \geq 0$.

It seems that the dependence among the variables of model \eqref{eq.KMP.1.1} has real impact on insurance practice. For this reason, the risk model \eqref{eq.KMP.1.1} attracted the interest of many researchers, and especially in case $d=2$. For bi-variate versions of  \eqref{eq.KMP.1.1} we refer to \cite{yang:li:2014}, \cite{jiang:wang:chen:xu:2015}, \cite{li:2024}, \cite{konstantinides:passalidis:2024a}, etc. For more dimensions, there only few references, see for example \cite{hult:lindskog:2006}, \cite{konstantinides:li:2016}, \cite{li:2016}, \cite{li:2022}, \cite{yang:su:2023}, \cite{cheng:konstantinides:wang:2024}, where is some of them the counting process is not necessarily common in the $d$ lines of business. In \cite{cheng:konstantinides:wang:2022} we find an extension of the time-dependent risk model from \cite{li:2016}, using the $\{\xi(t)\,,\;t\geq 0\}$, as a c\`{a}dl\`{a}g process, instead of L\'{e}vy process. Although, all these multivariate models provide a wide framework of dependence modeling, under the presence of heavy tails, there are three assumptions on the claims that can be regarded as restrictive. Namely
\begin{enumerate}

\item
The random vectors $\{{\bf X}^{(i)}\,,\;i\in \bbn\}$, which represent claims, belong to the class of multivariate regular variation, symbolically $MRV$. In spite of the importance of the $MRV$ class, in the frame of the the heavy-tailed distribution classes, it excludes several important distributions from the actuarial practice.

\item
Each claim-vector ${\bf X}^{(i)}$, has asymptotically dependent components. In this case, from one side it is restrictive to dependence, and from the other side it implies, in combination of the $MRV$, that the components of ${\bf X}^{(i)}$ follow distributions with tails that are strongly equivalent each other, which is clearly restrictive, especially for large values of $d$.

\item
The sequence of $\{{\bf X}^{(i)}\,,\;i\in \bbn\}$, contains independent and identically distributed terms.
\end{enumerate}

In this work we attempt to face the points $(1)$ - $(3)$. Most of the previous papers are focused on the points $(1)$ and $(2)$. Of course,  in order to avoid the points $(1)$ - $(3)$ we sacrifice the generality of the set family over the whole positive quadrant, except the point of origin, to some set from family $\mathscr{R}$, defined in Section $2$ below.

Concretely, we consider that the random vectors $\{{\bf X}^{(i)}\,,\;i\in \bbn\}$ are identically distributed, with multivariate distribution from class $(\mathcal{D} \cap \mathcal{L})_A$, larger that $MRV$, and we permit arbitrary dependence among the components of each vector ${\bf X}^{(i)}$. Further, we assume for the family of vectors $\{{\bf X}^{(i)}\,,\;i\in \bbn\}$ is subject to a weak dependence structure, which contains the independence as special case. This way, we obtain the interdependence among the claims, which means that each vector has dependent components (and here even arbitrary), and at the same time any two vectors are still under some dependence structure.

The restrictions in sets from the family $\mathscr{R}$, found in the seminal paper \cite{samorodnitsky:sun:2016}, can help to face the points $(1)$ and $(2)$, by defining the corresponding multivariate subexponentiality. Indeed, it was examined a multivariate renewal risk model, without interest rate, with independent, multivariate, subexponential claims. Next, in \cite{konstantinides:passalidis:2025g} the asymptotic behavior of the total discounted claims, in a Poisson risk model with multivariate subexponential claims was examined, under similar conditions, with the present ones, for the the price process of the investment portfolio. In last paper the time horizon is finite. However, in the two last works the point $(3)$ was not faced and furthermore, the asymptotic estimations are not uniform.

In this paper, we restrict the multivariate claim distribution class from the multivariate subexponentiality $\mathcal{S}_A$, into  the multivariate dominatedly varying, long tailed $(\mathcal{D} \cap \mathcal{L})_A$, in order to face point $(3)$, while the same time we establish uniformity in the asymptotic expressions.

The structure of the paper is as follows. In Section 2, we provide necessary preliminary statements for the distribution classes of multivariate heavy tailed distributions. In Section $3$ we provide the assumptions and the main result in finite time horizon, in which we obtain the asymptotic estimation of the behavior of the total, discounted claims of model \eqref{eq.KMP.1.1}, on a rare set $A$. As direct consequence of it, we can approximate the ruin probability on each rare set. In Section $4$, we generalize the results of Section $3$, in infinite time horizon. For this, we are forced to restrict ourselves in $MRV$ claim vectors, and we need some extra moment conditions for the stochastic process, that describes the investment portfolio prices. In spite of $MRV$, the points $(2)$ and $(3)$ are manageable. Finally, in Section $5$, we accommodate the proofs.

\section{Heavy-tailed random vectors} \label{sec.KLP.2}

For any set $A$ from the space $\bbr^d$, we denote by $A^c$ its complement set, by $\overline{A}$ its closed case and by $\partial A$ its border. For some event $E$, we denote by ${\bf 1}_E$ the indicator function of this set. For two positive numbers $a$ and $b$, we denote $a\vee b :=\max\{a,\,b\}$ and $a\wedge b :=\min\{a,\,b\}$.

All the vectors are $d$-dimensional, written in bold. For any two vectors, say ${\bf a}$ and ${\bf b}$, all the operations are defined component-wise. From now up, all the limit relations are understood as $\xto$, except it is said differently. For two $d$-variate, positive functions ${\bf f}$ and ${\bf g}$ and some set $A \in \bbr_+^d \setminus \{{\bf 0}\}=[0,\,\infty]^d \setminus \{{\bf 0}\}$, we denote ${\bf f}(x\,A)\sim c\,{\bf g}(x\,A)$, for some $c\in (0,\,\infty)$, or ${\bf f}(x\,A)= o[{\bf g}(x\,A)]$, or ${\bf f}(x\,A)= O[{\bf g}(x\,A)]$, if it holds
\beao
\lim \dfrac{{\bf f}(x\,A)}{{\bf g}(x\,A)}=c\,,\qquad \lim \dfrac{{\bf f}(x\,A)}{{\bf g}(x\,A)}=0\,,\qquad \limsup \dfrac{{\bf f}(x\,A)}{{\bf g}(x\,A)}< \infty\,,
\eeao
respectively. Furthermore, we write ${\bf f}(x\,A)\asymp {\bf g}(x\,A)$, if both relations ${\bf f}(x\,A)= O[{\bf g}(x\,A)]$ and  ${\bf g}(x\,A)= O[{\bf f}(x\,A)]$ hold. Correspondingly, for the $(d+1)$-variate, positive functions  ${\bf f}^*$ and ${\bf g}^*$, we denote ${\bf f}^*(x\,A;\,t)\sim {\bf g}^*(x\,A;\,t)$, uniformly for any $t \in \Delta$, for some non-empty set $\Delta$ if it holds
\beao
\lim\sup_{t \in \Delta} \left|\dfrac{{\bf f}^*(x\,A;\,t)}{{\bf g}^*(x\,A;\,t)} -1 \right| = 0\,.
\eeao  
Additionally, we write ${\bf f}^*(x\,A;\,t) \lesssim {\bf g}^*(x\,A;\,t)$, or ${\bf g}^*(x\,A;\,t) \gtrsim {\bf f}^*(x\,A;\,t)$, uniformly for any $t \in \Delta$, if
\beao
\limsup \sup_{t\in \Delta} \dfrac{{\bf f}^*(x\,A;\,t)}{{\bf g}^*(x\,A;\,t)} \leq 1\,.
\eeao
Finally, for any uni-variate distribution $V$, we denote its tail by $\overline{V}(x) = 1- V(x)$, for any $x\in \bbr$.

Now, we give the necessary preliminary concepts for the multivariate distributions with heavy tail, for later use. Let us notice, that a uni-variate distribution $V$ has heavy tail, if for any $\vep>0$ it holds
\beao
\int_0^{\infty} e^{\vep\,x}\,V(dx) = \infty\,.
\eeao
We have to remind that all random vectors are defined in the positive orthant of $\bbr_+^d\setminus \{{\bf 0}\}$, with each uni-variate distribution to have right endpoint equal to infinity, namely $\bV(x)>0$, for any $x \in \bbr$. It is worth to consider the several attempts to extend the subexponential property to multivariate distributions, as we know four different ones, by \cite{cline:resnick:1992}, \cite{omey:2006}, \cite{samorodnitsky:sun:2016}, \cite{konstantinides:passalidis:2024c}. In the case of last two, seems to operate complementary, while in the third one, we find a constructive criticism to the first two. 

In this work, we focus on the subexponential approach by \cite{samorodnitsky:sun:2016}, where it is based on the following fundamental set family
\beam \label{eq.KMP.2.1}
\mathscr{R}:=\left\{A \subsetneq \bbr^d\;:\; A \;\text{open,\,increasing}\,,\;A^c\; \text{convex}\,,\;{\bf 0} \notin \bA \right\}\,,
\eeam
where a set $A$ is named increasing, if for any ${\bf x} \in A$ and ${\bf y} \in \bbr_+^d$, it holds ${\bf x}+{\bf y}\in A$. By \cite[Lem. 4.5]{samorodnitsky:sun:2016}, we get that for any $A \in \mathscr{R}$, the random variable 
\beam \label{eq.KMP.2.2}
Y_A := \sup \{u\;:\;{\bf X} \in u\,A \}\,,
\eeam 
with ${\bf X}$ some random vector on $\bbr_+^d$, whose distribution is $F$, follows a proper distribution $F_A$, and its tail is given by
\beam \label{eq.KMP.2.3}
\bF_A(x) = \PP[{\bf X} \in x\,A ]= \PP\left[\sup_{{\bf p}\in I_A} {\bf p}^T\,{\bf X} >x\right]\,,
\eeam
for some index set $I_A \subsetneq \bbr^d$, see \cite[Lem. 4.3(c)]{samorodnitsky:sun:2016}, with ${\bf p}^T$ the inverse of vector ${\bf p}$. Thus, through relations \eqref{eq.KMP.2.1} - \eqref{eq.KMP.2.3}, is defined the multivariate subexponentiality on set $A$, symbolically $\mathcal{S}_A$.

We say that $F \in \mathcal{S}_A$, if $F_A \in \mathcal{S}$, namely for any (or equivalently, for some) integer $n\geq 2$ it holds
\beao
\lim \dfrac{\overline{F_A^{n*}}(x)}{\bF_A (x)}=n\,,
\eeao
where $F_A^{n*}$, represents the $n$-fold convolution of distribution $F_A$ with itself. Taking into account this definition of subexponentiality, in \cite{konstantinides:passalidis:2025g}, was introduced the following three distribution classes for some $A \in \mathscr{R}$. 

We say that distribution $F$ follows a multivariate long tail distribution on $A$, symbolically $F \in \mathcal{L}_A$, if $F_A\in \mathcal{L}$, namely for any (or equivalently, for some) $a>0$ it holds
\beao
\lim \dfrac{\overline{F_A}(x-a)}{\bF_A (x)}=1\,.
\eeao

We say that distribution $F$ follows a multivariate dominatedly varying distribution on $A$, symbolically $F \in \mathcal{D}_A$, if $F_A\in \mathcal{D}$, namely for any (or equivalently, for some) $b \in (0,\,1)$ it holds
\beao
\limsup \dfrac{\overline{F_A}(b\,x)}{\bF_A (x)} < \infty\,.
\eeao
Furthermore is defined the property $F \in (\mathcal{D} \cap \mathcal{L})_A$, if $F_A \in \mathcal{D} \cap \mathcal{L}$. We note that, from \cite{goldie:1978} we find that $\mathcal{D} \cap \mathcal{L} = \mathcal{D} \cap \mathcal{S}$. This way, for any distribution class $\mathcal{B} \in \{\mathcal{S},\,\mathcal{L},\,\mathcal{D},\,\mathcal{D} \cap \mathcal{L}\}$, we denote $\mathcal{B}_{\mathscr{R}}:=\bigcap_{A\in \mathscr{R}} \mathcal{B}_A$. By the definition of these classes, is directly implied that remain invariant the ordering of the uni-variate classes to the multivariate ones.

Next, we need the regular variation for uni-variate distributions. We say that a distribution $V$, belongs to the class of regularly varying distributions with some index $\alpha \in (0,\,\infty)$, symbolically $V \in \mathcal{R}_{-\alpha}$, if for any $t>0$ it holds
\beao
\lim \dfrac{\bV(t\,x)}{\bV(x)} = t^{-\alpha}\,.
\eeao

The class of multivariate regular varying distributions, symbolically $MRV$, was introduced in \cite{dehaan:resnick:1981}, and represents the most famous multivariate distributions with heavy tail. We provide its standard form, which is rather flexible and common, in comparison to the non-standard one. We say that a random vector ${\bf X}$ with distribution $F$, belongs to standard $MRV$ class, if there exists a distribution $V \in \mathcal{R}_{-\alpha}$, with $\alpha \in (0,\,\infty)$, and a non-degenerated to zero Radon measure $\mu$, such that it holds
\beam \label{eq.KMP.2.8}
\dfrac{1}{\bV(x)}\,\PP[{\bf X} \in x\,\bbb] \stackrel{v}{\rightarrow} \mu(\bbb)\,,
\eeam
as $\xto$, for any $\mu$-continuous Borel set $\bbb \in \bbr_+^d\setminus \{{\bf 0}\} $, where the "$\stackrel{v}{\rightarrow} $" denotes the vague convergence. For distributions satisfying \eqref{eq.KMP.2.8}, we write $F \in MRV(\alpha,\,V,\,\mu)$. Further we know that the measure $\mu$, is positively homogeneous, namely for any Borel set $\bbb \in \bbr_+^d\setminus \{{\bf 0}\} $ and any $c>0$, it holds
\beao
\mu\left(c^{1/\alpha}\,\bbb \right) = \dfrac 1c\,\mu(\bbb)\,.
\eeao
We refer to \cite{resnick:2007}, \cite{mikosch:wintenberger:2024}, \cite{resnick:2024} for further treatments in $MRV$. 

In relation to points $(1)$ - $(3)$ in first section, we remark that the asymptotic dependence among the components of ${\bf X}$, is understood in the case of $MRV$ as follows
\beao
\mu\left((1,\,\infty]\times \cdots \times (1,\,\infty] \right) > 0\,.
\eeao
Next, it is easy to see, that $F \in MRV(\alpha,\,V,\,\mu)$, with $\alpha \in (0,\,\infty)$, then for any $A \in \mathscr{R} $, it holds $F_A \in \mathcal{R}_{-\alpha} $, with the same regular variation index, but the inverse implication does not necessarily hold. Even more, according to \cite[Prop. 4.14]{samorodnitsky:sun:2016} and \cite[Prop. 2.1]{konstantinides:passalidis:2025g}, it holds 
\beao
\bigcup_{0< a < \infty} MRV(\alpha,\,V,\,\mu) \subsetneq (\mathcal{D}\cap \mathcal{L})_{\mathscr{R}}\subsetneq \mathcal{S}_{\mathscr{R}} \subsetneq \mathcal{L}_{\mathscr{R}}\,.
\eeao

In order to formulate the main result, we provide a classical kind of index for uni-variate distributions, with heavy tail. The upper Matuszewska index, for a distribution $V$, with infinite right endpoint, is given as
\beam \label{eq.KMP.2.12}
J_{V}^+:=- \lim_{v\to \infty} \dfrac {\ln \bV_{*}(v)}{\ln v}\,, 
\eeam
with
\beao
\bV_{*}(v):=\liminf_{\xto} \dfrac {\bV(v\,x)}{\bV(x)}\,,
\eeao
for any $v>1$, thus from \eqref{eq.KMP.2.12} we see that $0\leq J_V^+ \leq \infty$, for any distribution $V$ with right endpoint. It is well-known that $V \in \mathcal{D}$ is equivalent to $J_{V}^+ < \infty$, and further if $V \in \mathcal{D}$ then for any $p>J_{V}^+$, it holds 
\beam \label{eq.KMP.2.13}
x^{-p}=o[\bV(x)]\,, 
\eeam
see for example \cite[Sec. 2.4]{leipus:siaulys:konstantinides:2023} or \cite[Lem.3.5]{tang:tsitsiashvili:2003}. Let us observe that if $V \in \mathcal{R}_{-\alpha}$, for some $\alpha \in (0,\,\infty)$, then $J_V^+ = \alpha$. For more distribution properties with heavy tails, in uni-variate case, see in \cite{foss:korshunov:zachary:2013}, \cite{konstantinides:2018}, etc.

\section{Uniform asymptotic estimates over finite intervals} \label{subsec.KMP.2.2}

Before formulation of assumptions for the main results, we should remark that all the results established for some set $A\in \mathscr{R}$, hold also for the set family $\mathscr{R}$. Further, we define the set $\Lambda = \{t\;:\;\lambda(t) >0\}$ and $\underline{t}:= \inf \{t\;:\;\lambda(t) >0\}$. Hence, we obtain
\beao
\Lambda= 
\begin{cases}
[\underline{t},\,\infty]  & \text{if } \;\PP[T_1=\underline{t}]>0\,,\\[2mm]
(\underline{t},\,\infty]\,, & \text{if } \;\PP[T_1=\underline{t}]=0\,.
\end{cases} 
\eeao
Now we denote, for some fixed $T \in (0,\,\infty)$, $\Lambda_T=(0,\,T]\cap \Lambda$. Finally, we represent the vector of the discounted aggregate claims, up to the moment $t\geq 0$, as 
\beao
{\bf D}(t)= \sum_{i=1}^{N(t)} {\bf X}^{(i)}\,e^{-\xi(T_i)}\,,
\eeao
and we look for uniform asymptotic estimations of the probability $\PP[{\bf D}(t) \in x\,A]$, as $\xto$. The set $x\,A$ can be understood as a rare-set, which can take various forms, that can have sense in risk theory, see for example in \cite[Sec. 4,\,5]{samorodnitsky:sun:2016} and Remark \ref{rem.KMP.3.1} below. Also for the case of the ruin set $L$, see in Remark \ref{rem.KMP.3.2} below, that can be seen as related to the sets $A$.

Now, we proceed to the first assumption for the risk model \eqref{eq.KMP.1.1}, that is related to the yield processes of the investment portfolio.

\begin{assumption} \label{ass.KMP.2.1}
Let there exist constants $C_2 \geq C_1 \geq 0$, which are dependent on the $T \in (0,\,\infty)$, such that it holds
\beao
\PP\left[\inf_{0\leq t \leq T} \xi(t) \geq -C_1\right] =1\,,\qquad \qquad \PP\left[\sup_{0\leq t \leq T} \xi(t) \leq C_2\right] =1\,,
\eeao
\end{assumption}
Due to finite time horizon, Assumption \ref{ass.KMP.2.1} is satisfied by c\'{a}dl\'{a}g processes, that means continuous from right with finite limits from left, which include the L\`{e}vy processes as a special case. 

Second assumption has to do with the dependence among claim vectors, which are with arbitrarily dependent components. This dependence is based on the tail asymptotic independence of $\{Y_A^{(i)}\,,\;i \in \bbn\}$, introduced in \cite{geluk:tang:2009} . In present paper, we say that the random vectors ${\bf X}^{(i)}$ and ${\bf X}^{(j)}$ are independent on $A\in \mathscr{R}$, if the random variables $Y_A^{(i)} := \sup \{u\;:\;{\bf X}^{(i)} \in u\,A \}$ and $Y_A^{(j)} := \sup \{u\;:\;{\bf X}^{(j)} \in u\,A \}$ are independent.

\begin{assumption} \label{ass.KMP.2.2}
Let $A\in \mathscr{R}$. The sequence of claim vectors $\{{\bf X}^{(i)}\,,\;i\in \bbn\}$ has identically distributed terms, with common distribution $F \in (\mathcal{D}\cap \mathcal{L})_A$. Furthermore, the ${\bf X}^{(i)}$ and ${\bf X}^{(j)}$, for $i\neq j$, are tail asymptotically independent on $A$, symbolically $TAI_A$, in the sense that the $Y_A^{(i)}$ and $Y_A^{(j)}$ are tail asymptotically independent, namely
\beam \label{eq.KMP.2.17}
\lim_{x_i\wedge x_j \to \infty} \PP\left[Y_A^{(i)} > x_i\;|\;Y_A^{(j)} > x_j\right] =0\,,
\eeam
with $x_i \wedge x_j := \min\{x_i,\,x_j\}$.
\end{assumption}
We can observe that relation \eqref{eq.KMP.2.17}, contains as special case the independence, hence also the independence on $A$ for the vectors ${\bf X}^{(i)}$ and ${\bf X}^{(j)}$. 

Assumption \ref{ass.KMP.2.3} indicates the independence among the premiums, the financial risks, the claim sizes and their counting process.

\begin{assumption} \label{ass.KMP.2.3}
The $\{{\bf X}^{(i)}\,,\;i\in \bbn \}$,  $\{N(t)\,,\;t\geq 0 \}$, $\{\xi(t)\,,\;t\geq 0 \}$ and $\{{\bf c}(t)\,,\;t\geq 0 \}$ are mutually independent.
\end{assumption}

We are in position to give the first main result. We should notice, that in case $d=1$, when $A = (1,\,\infty)$, relation \eqref{eq.KMP.2.18}, generalizes \cite[Th. 2.1]{yang:wang:konstantinides:2014}, with respect to price processes of the investment portfolio.

\bth \label{th.KMP.2.1}
Let $A\in \mathscr{R}$. We consider the risk model \eqref{eq.KMP.1.1} and suppose that Assumptions \ref{ass.KMP.2.1}, \ref{ass.KMP.2.2} and \ref{ass.KMP.2.3} are true. Then for any fixed $T>0$, it holds
\beam \label{eq.KMP.2.18}
 \PP\left[{\bf D}(t) \in x\,A \right] \sim \int_0^t \PP[{\bf X}\,e^{-\xi(s)} \in x\, A]\,\lambda(ds)\,,
\eeam
uniformly for any $t \in \Lambda_T$.
\ethe

\bre \label{rem.KMP.3.1}
Before proceeding to ruin probability, we state that two important forms of $A\in \mathscr{R}$, in the actuarial practice are as follow
\beam \label{eq.KMP.2.19}
A=\left\{ {\bf x} \;:\; x_i> b_i\,,\;\exists i=1,\,\ldots,\,d \right\}\,,
\eeam
where $b_1,\,\ldots,\,b_d >0$ and 
\beam \label{eq.KMP.2.20}
A=\left\{ {\bf x} \;:\; \sum_{i=1}^d l_i\,x_i > u \right\}\,,
\eeam
where $u>0$ and 
\beao
\sum_{i=1}^d l_i =1\,,
\eeao 
with $l_1,\,\ldots,\,l_d \geq 0$.
\ere
Therefore, if the set $A$ takes the form of \eqref{eq.KMP.2.19}, then \eqref{eq.KMP.2.18} is interpreted as the asymptotic probability of the event the discounted aggregate claims to have a line of business, which would have larger claim than its  initial capital, while if the set $A$ takes the form of \eqref{eq.KMP.2.20}, then \eqref{eq.KMP.2.18} is interpreted as the asymptotic probability of the event the discounted aggregate claims to have a total sum, which would exceed the total initial capital.

The next assumption presents some basic characteristics of the ruin sets, that we need. Let us note that the set $L$ is called decreasing, if the $-L$ is increasing.

\begin{assumption} \label{ass.KMP.2.4}
Let $L$ a ruin set, which is open, decreasing, with $L^c$ convex and ${\bf 0} \in \partial L$. Additionally, we suppose that for any $x>0$, it holds $L = x\,L$.
\end{assumption}  

\bre \label{rem.KMP.3.2}
Assumption \ref{ass.KMP.2.4} can be found in \cite{hult:lindskog:2006} and used in \cite{samorodnitsky:sun:2016}. Following the line of discussion in \cite[Sec. 5]{samorodnitsky:sun:2016}, the set $A=({\bf l} - L) \in \mathscr{R}$, and further two improtant cases of ruin sets, that satisfy Assumption \ref{ass.KMP.2.4}, are as follow
\beao
L=\left\{ {\bf x} \;:\; x_i< 0\,,\;\exists i=1,\,\ldots,\,d \right\}\,,
\eeao
where the entrance of the surplus in this ruin set, means that bankrupt in one of the $d$ lines of business and 
\beam \label{eq.KMP.2.21b}
L=\left\{ {\bf x} \;:\; \sum_{i=1}^d x_i < 0 \right\}\,,
\eeam
where the entrance of the surplus in this ruin set, means that bankrupt of all the $d$ lines of business in total. We should also notice that for $d=1$, the $L=(-\infty,\,0)$ represents a ruin set that satisfies Assumption \ref{ass.KMP.2.4}, see in \cite{hult:lindskog:2006} for more information about ruin sets, satisfying  Assumption \ref{ass.KMP.2.4}.
\ere

Next, we introduce the ruin probability over finite time horizon, with respect to entrance probability of the surplus into some ruin set, up to time $t>0$, and through the set 
\beao
A=({\bf l} - L) \in \mathscr{R}\,.
\eeao 
We obtain
\beao
&&\psi_{{\bf l},L}(x;\,t)\\[2mm] 
&&:=\PP[{\bf U}(s) \in L,\;\exists \;s\in (0,\,t]] =\PP\left[{\bf D}(s) -\int_0^s e^{-\xi(y)}{\bf c}(y)dy \in x\,({\bf l}-L),\;\exists s \in (0,\,t] \right]\\[2mm]
&&=\PP\left[{\bf D}(s) -\int_0^s e^{-\xi(y)}\,{\bf c}(y)\,dy \in x\,A\,,\;\exists s \in (0,\,t] \right]\,,
\eeao
where 
\beao
\int_0^t e^{-\xi(y)}{\bf c}(y)dy :=
\left( 
\begin{array}{c}
\int_{0-}^{t} e^{-\xi(y)}\,c_1(y)\,dy \\ 
\vdots \\ 
\int_{0-}^{t} e^{-\xi(y)}\,c_d(y)\,dy 
\end{array} 
\right) \,,
\eeao
for all $t>0$. Now we are in position to estimate the ruin probability over finite time horizon in the risk model \eqref{eq.KMP.1.1}.

\bco \label{cor.KMP.2.1}
Let $A=({\bf l} - L) \in \mathscr{R}$. Under the conditions of Theorem \ref{th.KMP.2.1}, for the risk model \eqref{eq.KMP.1.1}, for any fixed $T>0$ it holds
\beam \label{eq.KMP.2.24} 
\psi_{{\bf l},L}(x;\,t)\sim \int_0^t \PP\left[{\bf X}\,e^{-\xi(s)} \in x\,A \right]\,\lambda(ds)\,,
\eeam
uniformly for any $t \in \Lambda_T$.
\eco

\section{Global uniformity}

In this section we study the risk model \eqref{eq.KMP.1.1} and provide uniformly asymptotic estimations over the whole $\Lambda$, for the asymptotic behavior of the discounted aggregate claims and also the ruin probability. In relation to Section 3, that is restricted to finite time horizon, here through some harder conditions we achieve the estimations over the whole time horizon. The first of conditions is the restriction of the distribution class to $MRV$ from the $(\mathcal{D} \cap \mathcal{L})_A$, while in the second condition we adopt Assuption \ref{ass.KMP.4.1} below, instead of Assumption \ref{ass.KMP.2.1}. The $\alpha \in (0,\,\infty)$ denotes the regular variation index of the random vector ${\bf X}$.

\begin{assumption} \label{ass.KMP.4.1}
Let $\{\xi(t)\,,\;t\geq 0 \}$ a c\'{a}dl\'{a}g process, satisfying the conditions
\begin{enumerate}
\item
If $\alpha \in [1,\,\infty)$, then there exist some $k_1,\,k_2$, with $0<k_1<\alpha < k_2 < \infty$, such that
\beam \label{eq.KMP.4.1} 
\sum_{i=1}^{\infty} \E\left[\left(e^{-\xi(T_i)} \right)^{k_1} \right]^{1/k_1}\vee \E\left[\left(e^{-\xi(T_i)} \right)^{k_2} \right]^{1/k_2}< \infty\,.
\eeam

\item
If $\alpha \in (0,\,1)$, then there exist some $k_1^*,\,k_2^*$, with $0<k_1^*<\alpha < k_2^* < 1$, such that
\beam \label{eq.KMP.4.2} 
\sum_{i=1}^{\infty} \E\left[\left(e^{-\xi(T_i)} \right)^{k_1^*} \right] \vee \E\left[\left(e^{-\xi(T_i)} \right)^{k_2^*} \right]< \infty\,.
\eeam
\end{enumerate}
\end{assumption} 

Assumption \ref{ass.KMP.4.1}, in fact contains a large variety of processes, which describe the return process of the investment portfolio. In what follows, we present an example, that satisfies Assumption \ref{ass.KMP.4.1}. Indeed, this assumption holds for some geometric L\`{e}vy processes. 

\bexam \label{exam.KMP.4.1}
Let $\{\xi(t)\,,\;t\geq 0\}$ be a L\`{e}vy process, with Laplace exponent
\beao
\phi(z) = \ln \,\E\left[e^{-z\,\xi(1)}\right]\,,
\eeao
for any $z \in \bbr$. Thus, when the Laplace exponent is finite $\phi(z)<\infty$, we obtain 
\beao
\E\left[e^{-z\,\xi(t)} \right]=e^{t\,\phi(s)} \,,
\eeao
for any $t \in \bbr$, see \cite[Ch. 3]{cont:tankov:2004}. From the fact that the Laplace exponent $\phi(z)$ is convex function with respect to $z$, and also that $\phi(0)=0$, we find that if for some $k_2^*>0$, it holds $\phi(k_2^*)<0$, then for any $x \in (0,\,k_2^*]$ we get $\phi(x)<0$.

Now, we show relation \eqref{eq.KMP.4.2} in Assumption \ref{ass.KMP.4.1}, with $0<k_1^* < \alpha <k_2^* <1$. Let $\phi(k_2^*)<0$, that means also $\phi(k_1^*)<0$, then since $\{N(t)\,,\;t>0\}$ is a renewal process, is implied 
\beao
\E\left[e^{-k_j^*\,\xi(T_i)}\right] = \left( \E\left[e^{\phi(k_j^*)\,T_1} \right]\right)^i\,,
\eeao
for $j=1,\,2$. From the previous statements we obtain
\beao
&&\sum_{i=1}^{\infty} \E\left[\left(e^{-\xi(T_i)} \right)^{k_1^*} \right]\vee \E\left[\left(e^{-\xi(T_i)} \right)^{k_2^*} \right]\leq \sum_{i=1}^{\infty} \left(\E\left[ e^{-k_1^*\,\xi(T_i)}  \right]+ \E\left[e^{-k_2^*\,\xi(T_i)} \right] \right)\\[2mm]
&&\leq \sum_{i=1}^{\infty} \left(\E\left[ e^{\phi(k_1^*)\,T_i} \right]+ \E\left[e^{\phi(k_2^*)\,T_i} \right] \right) < \infty\,,
\eeao
which verifies relation \eqref{eq.KMP.4.2}. Relation \eqref{eq.KMP.4.1} is satisfied if we assume that $\phi(k_2) <0$, by following similar steps, hence we omit the argument. 
\eexam

We are now in position to provide the main result of this section. We observe that except the asymptotic estimation for infinite time horizon, we also obtain a more direct asymptotic expression in comparison to \eqref{eq.KMP.2.24}, which is due to a kind of multivariate version of the Breiman's lemma, see further in Subsection 5.2.

\bth \label{th.KMP.4.1}
Let $A \in \mathscr{R}$ and hold Assumption \ref{ass.KMP.4.1}, together with Assumption \ref{ass.KMP.2.2} and Assumption \ref{ass.KMP.2.3}, under the restriction that the claim distribution $F$ belongs to $MRV$,
\beao
F \in MRV(\alpha,\,V,\,\mu)\,,
\eeao
with $\alpha \in (0,\,\infty)$. Then, in risk model \eqref{eq.KMP.1.1} the relation 
\beam \label{eq.KMP.4.3} 
\PP\left[{\bf D}(t) \in x\,A \right] \sim \PP\left[{\bf X} \in x\,A  \right]\,\int_0^t \E\left[ e^{-\alpha\,\xi(s)}\right]\,\lambda(s)\,.
\eeam
holds uniformaly for any $t \in \Lambda$.
\ethe

We should mention that in Theorem \ref{th.KMP.4.1}, as well as in Corollary \ref{cor.KMP.4.1}, we do not assume a dependence structure for the components of each claim vector, and we keep  the $TAI_A$ dependence structure among the claims, thence our results do not face point $(1)$, but face points $(2)$ and $(3)$, from Section 1, although we are forced to be restricted only in set from family $\mathscr{R}$.
	
In the following consequence, we provide the ruin probability over infinite time horizon.

\bco \label{cor.KMP.4.1}
Let $A=({\bf l} - L) \in \mathscr{R}$. Under the conditions of Theorem \ref{th.KMP.4.1} and under the assumption that 
\beao
\int_0^t  e^{-\xi(y)}\,c_{i}(y)\,dy <\infty\,,
\eeao
for any $t \in \Lambda$, and any $i=1,\,\ldots,\,d$, then for the risk model \eqref{eq.KMP.1.1} it holds
\beam \label{eq.KMP.4.5} 
\psi_{{\bf l},L}(x;\,t)\sim \PP\left[{\bf X} \in x\,A  \right]\,\int_0^t \E\left[e^{-\alpha\,\xi(s)} \right]\,\lambda(ds)\,,
\eeam
uniformly for any $t \in \Lambda$.
\eco

\section{Argumentation}

\subsection{Finite time horizon}
Let start now the proofs.

{\bf Proof of Theorem \ref{th.KMP.2.1}.}~
At first from Assumption \ref{ass.KMP.2.1}, is implied that for any $t \in \Lambda_T$, it holds $0 < e^{-C_2} \leq e^{-\xi(t)} \leq e^{C_1} < \infty$. Now, for any integer $m \in \bbn$ and any $x>0$, it holds
\beam \label{eq.KMP.3.1}
&&\PP\left[{\bf D}(t) \in x\,A \right]=\PP\left[\sum_{i=1}^{N(t)} {\bf X}^{(i)}\,e^{-\xi(T_i)} \in x\,A  \right]\\[2mm] \notag
&&=\left( \sum_{n=1}^{m} + \sum_{n=m+1}^{\infty} \right)\PP\left[\sum_{i=1}^{n} {\bf X}^{(i)}\,e^{-\xi(T_i)} \in x\,A\,,\;N(t)=n \right]=J_1(x,\,m,\,t)+J_2(x,\,m,\,t)\,.
\eeam 

Let calculate first the $J_2(x,\,m,\,t)$. Using \cite[Prop. 2.4]{konstantinides:passalidis:2025g}, via Markov's inequality, choosing some $p>\alpha_{F_A}$, for $x > m$, we obtain
\beam \label{eq.KMP.3.2} \notag
&&J_2(x,\,m,\,t) \leq \sum_{n=m+1}^{\infty} \PP\left[\sum_{i=1}^{n} Y_A^{(i)}\,e^{-\xi(T_i)} > x\,,\;N(t)=n \right]\\[2mm]
&&=\left( \sum_{m<n\leq x} + \sum_{n>x}^{\infty} \right)\PP\left[\sum_{i=1}^{n} Y_A^{(i)}\,e^{-\xi(T_i)} > x \right]\,\PP[N(t)=n]\\[2mm] \notag
&&\leq\sum_{m<n\leq x} \PP\left[\sum_{i=1}^{n} Y_A^{(i)}\,e^{C_1} > x \right]\,\PP[N(t)=n]+\PP[N(t) > x]\\[2mm] \notag
&&\leq\sum_{m<n\leq x} n\,\PP\left[Y_A^{(i)}\,e^{C_1} > \dfrac xn \right]\,\PP[N(t)=n]+x^{-(p+1)}\,\E\left[N^{p+1}(t)\,{\bf 1}_{\{N(t)> x\}}\right]\\[2mm] \notag
&&\leq K_1 \PP[Y_A>x]\,\sum_{m<n\leq x} n^{p+1}\,\PP[N(t)=n]+x^{-(p+1)}\,\E\left[N^{p+1}(t)\,{\bf 1}_{\{N(t)> x\}}\right]\\[2mm] \notag
&&\lesssim K_1 \PP[{\bf X} \in x\,A]\,\E\left[N^{p+1}(t)\,{\bf 1}_{\{N(t)> m\}}\right]\,,
\eeam 
with $Y_A^{(i)}=\sup\{u\;:\;{\bf X}^{(i)} \in u\,A\}$, where in the fifth step the constant $K_1>0$ follows from the membership of the distribution of $Y_A$ to class $\mathcal{D}$, while in the last step we used \eqref{eq.KMP.2.13}. Next, due to the class $\mathcal{D}$ properties, we find for any $t \in \Lambda_T$ the relation
\beam \label{eq.KMP.3.3}
\int_0^t \PP\left[ {\bf X} \,e^{-\xi(s)} \in x\,A \right]\,\lambda(ds) &=& \int_0^t \PP\left[ Y_A\,e^{-\xi(s)}> x \right]\,\lambda(ds) \\[2mm] \notag
&\geq& \PP\left[ Y_A\,e^{-C_2} > x \right]\,\lambda(t) \asymp \PP\left[ {\bf X} \in  x\,A\right]\,\lambda(t) \,.
\eeam
By \eqref{eq.KMP.3.3} we conclude that there exists some constant $K_2>0$, such that it holds
\beam \label{eq.KMP.3.4}
\int_0^t \PP\left[ {\bf X} \,e^{-\xi(s)} \in x\,A \right]\,\lambda(ds) \geq K_2\, \PP\left[ {\bf X} \in  x\,A\right]\,\lambda(t) \,,
\eeam
for sufficient large $x>0$ and for any $t \in \Lambda_T$. Hence, by relations \eqref{eq.KMP.3.2} and \eqref{eq.KMP.3.4}, we obtain
\beam \label{eq.KMP.3.5} \notag
&&\lim_{m\to \infty} \limsup_{\xto} \sup_{t \in \Lambda_T} \dfrac {J_2(x,\,m,\,t)}{\int_0^t\PP\left[ {\bf X} \,e^{-\xi(s)} \in x\,A \right]\,\lambda(ds)} \\[2mm] \notag
&&\leq \lim_{m\to \infty} \limsup_{\xto} \sup_{t \in \Lambda_T} \dfrac {K_1 \PP[{\bf X} \in x\,A]\,\E\left[N^{p+1}(t)\,{\bf 1}_{\{N(t)> m\}}\right])}{ K_2\, \PP\left[ {\bf X} \in  x\,A\right]\,\lambda(t)} \\[2mm]
&&= \lim_{m\to \infty} \sup_{t \in \Lambda_T}  \dfrac {K_1}{K_2\,\lambda(t)}\,\E\left[N^{p+1}(t)\,{\bf 1}_{\{N(t)> m\}}\right]=0\,,
\eeam
where in last step we took advantage of \cite[Lem. 3.2]{tang:2007}.

Now we estimate the $J_1(x,\,m,\,t)$. At first, we obtain that for $1\leq n \leq m$, for any $t \in \Lambda_T$ it holds
\beam \label{eq.KMP.3.6}
&& \PP\left[\sum_{i=1}^{n}  {\bf X}^{(i)}\,e^{-\xi(T_i)} \in x\,A\,,\;N(t)=n \right]\\[2mm] \notag
&&=\int_{\{0\leq t_1 \leq \cdots \leq t_n \leq t,\, t_{n+1}>t\}} \PP\left[ \sum_{i=1}^{n}  {\bf X}^{(i)}\,e^{-\xi(t_i)} \in x\,A \right]\,\PP\left[ T_1\in dt_1,\,\ldots,\,T_{n+1}\in dt_{n+1} \right] \,.
\eeam
From \cite[Prop. 2.4]{konstantinides:passalidis:2025g} together with \cite[Th. 2.1]{li:2013} we find
\beao
\PP\left[\sum_{i=1}^{n}  {\bf X}^{(i)}\,e^{-\xi(t_i)} \in x\,A \right]&\leq&  \PP\left[\sum_{i=1}^{n} Y_A^{(i)}\,e^{-\xi(t_i)} > x \right] \\[2mm] \notag
&\sim& \sum_{i=1}^{n} \PP\left[Y_A^{(i)}\,e^{-\xi(t_i)} > x \right]=\sum_{i=1}^{n} \PP\left[{\bf X}^{(i)}\,e^{-\xi(t_i)} \in x\,A \right] \,,
\eeao
uniformly for any  $t_i \in \Lambda_T$, where the pre-last step permits the application of \cite[Th.2.1]{li:2013}, because of Assumption \ref{ass.KMP.2.1} and the $TAI_A$ property. Hence, for some $\delta_1>0$, there exists some $x_0$, large enough, independent of $t$, such that, for any $x\geq x_0$, and for $0\leq t_1 \leq \cdots \leq t_n \leq t$ and $1\leq n \leq m$, it holds
\beam \label{eq.KMP.3.8}
\PP\left[\sum_{i=1}^{n}  {\bf X}^{(i)}\,e^{-\xi(t_i)} \in x\,A \right]\leq (1+\delta_1) \sum_{i=1}^{n} \PP\left[{\bf X}^{(i)}\,e^{-\xi(t_i)} \in x\,A \right] \,.
\eeam

For the corresponding lower bound of the scale mixture sum in the left hand side of \eqref{eq.KMP.3.8}, due to $F_A \in \mathcal{D}$, by \cite[Th. 3.3(iv)]{cline:samorodnitsky:1994}, follows the relation
\beam \label{eq.KMP.3.9}
\PP\left[ {\bf X}^{(i)}\,e^{-\xi(t_i)} \in x\,A \right]=\PP\left[Y_A^{(i)}\,e^{-\xi(t_i)} > x \right]\asymp \PP\left[ Y_A^{(i)} > x \right]= \PP\left[{\bf X}^{(i)} \in x\,A \right] \,,
\eeam
and thus by relation \eqref{eq.KMP.3.9}, we find that there exist constants $0 < k_1(i) \leq k_2(i) $ and some large enough $x_0^*>0$, independent of $t$, such that, for any $x \geq x_0^*$ and $0\leq t_1 \leq \cdots \leq t_n \leq t$, it holds
\beam \label{eq.KMP.3.10}
k_1(i) \leq \dfrac{\PP\left[ {\bf X}^{(i)}\,e^{-\xi(t_i)} \in x\,A \right]}{\PP\left[ {\bf X}^{(i)} \in x\,A \right]}\leq k_2(i)\,.
\eeam

Now, because of $TAI_A$ dependence structure, for some $\vep>0$, there exists some large enough $\tilde{x}_0>0$, independent of $t$, such that for any $x \geq \tilde{x}_0$, all $i \neq j$ and all $0\leq t_1 \leq \cdots \leq t_n \leq t$, it holds 
\beam \label{eq.KMP.3.11}
&&\PP\left[ {\bf X}^{(i)}\,e^{-\xi(t_i)} \in x\,A\,,\; {\bf X}^{(j)}\,e^{-\xi(t_j)} \in x\,A\right] \\[2mm] \notag
&& \leq\PP\left[ {\bf X}^{(i)}\,e^{C_1} \in x\,A\,,\; {\bf X}^{(j)}\,e^{C_1} \in x\,A\right] \leq \vep\,\PP\left[{\bf X}^{(i)}\,e^{C_1} \in x\,A \right] \leq \vep\,v_1 \,\PP\left[{\bf X}^{(i)} \in x\,A \right]\,,
\eeam
where the constant $v_1>0$ follows by the fact that $F \in \mathcal{D}_A$. Since the ${\bf X}^{(i)}$ and ${\bf X}^{(j)}$ as also the $e^{-\xi(t)}$, for $t \in [0,\,T]$ are non-negative (and $A$ is increasing set), applying Bonferroni's inequality and using relations \eqref{eq.KMP.3.10} and \eqref{eq.KMP.3.11}, for any $x\geq x_0^* \vee \tilde{x}_0$, we obtain
\beam \label{eq.KMP.3.12} \notag
&&\PP\left[ \sum_{i=1}^{n}  {\bf X}^{(i)}\,e^{-\xi(t_i)} \in x\,A  \right]\geq \PP\left[ \cup_{i=1}^{n} \left\{ {\bf X}^{(i)}\,e^{-\xi(t_i)} \in x\,A \right\} \right] \\[2mm] \notag
&&\geq \sum_{i=1}^{n} \PP\left[  {\bf X}^{(i)}\,e^{-\xi(t_i)} \in x\,A  \right]- \sum_{i=1}^{n} \sum_{i<j\leq n}\PP\left[  {\bf X}^{(i)}\,e^{-\xi(t_i)} \in x\,A\,,\; {\bf X}^{(j)}\,e^{-\xi(t_j)} \in x\,A  \right]\\[2mm] \notag
&&\geq \sum_{i=1}^{n} \PP\left[  {\bf X}^{(i)}\,e^{-\xi(t_i)} \in x\,A  \right]- \vep\,v_1\,\sum_{i=1}^{n}\PP\left[  {\bf X}^{(i)}\in x\,A  \right]\geq \sum_{i=1}^{n} \PP\left[  {\bf X}^{(i)}\,e^{-\xi(t_i)} \in x\,A  \right]\\[2mm]
&&- \vep\,\sum_{i=1}^{n}\,\dfrac {v_1}{k_1(i)}\PP\left[ {\bf X}^{(i)}\,e^{-\xi(t_i)}\in x\,A  \right]\geq(1-u_1)\,\sum_{i=1}^{n}\PP\left[ {\bf X}^{(i)}\,e^{-\xi(t_i)}\in x\,A  \right]\,,
\eeam
with 
\beao
u_1=\vep\,\dfrac{v_1}{\wedge_{i=1}^{n} k_1(i)} >0\,,
\eeao 
that can be arbitrarily close to zero. Therefore, by relations \eqref{eq.KMP.3.8} and \eqref{eq.KMP.3.12}, for 
\beao
\delta:=\delta_1 \vee u_1> 0\,,
\eeao 
and for any $x\geq x_0 \vee x_0^* \vee \tilde{x}_0$, we obtain
\beam \label{eq.KMP.3.13} \notag
(1-\delta) \sum_{i=1}^{n} \PP\left[  {\bf X}^{(i)}\,e^{-\xi(t_i)} \in x\,A  \right]&\leq& \PP\left[ \sum_{i=1}^{n}  {\bf X}^{(i)}\,e^{-\xi(t_i)} \in x\,A  \right] \\[2mm]
&\leq& (1+\delta)\,\sum_{i=1}^{n}\PP\left[ {\bf X}^{(i)}\,e^{-\xi(t_i)}\in x\,A  \right]\,.
\eeam
Thus, by relations \eqref{eq.KMP.3.6} and \eqref{eq.KMP.3.13} we conclude
\beao
&&(1-\delta) \sum_{i=1}^{n} \PP\left[  {\bf X}^{(i)}\,e^{-\xi(T_i)} \in x\,A \,,\;N(t)=n \right]\leq \PP\left[ \sum_{i=1}^{n}  {\bf X}^{(i)}\,e^{-\xi(T_i)} \in x\,A\,,\;N(t)=n  \right] \\[2mm]
&&\leq (1+\delta)\,\sum_{i=1}^{n}\PP\left[ {\bf X}^{(i)}\,e^{-\xi(T_i)}\in x\,A\,,\;N(t)=n  \right]\,.
\eeao
Hence, letting $\delta$ to tend to zero, as $\xto$, we find the asymptotic equivalence
\beam \label{eq.KMP.3.16} \notag
\PP\left[ \sum_{i=1}^{n}  {\bf X}^{(i)}\,e^{-\xi(T_i)} \in x\,A\,,\;N(t)=n  \right]\sim \sum_{i=1}^{n} \PP\left[  {\bf X}^{(i)}\,e^{-\xi(T_i)} \in x\,A \,,\;N(t)=n \right]\,,
\eeam
uniformly for any $t \in \Lambda_T$. 

Now, from \eqref{eq.KMP.3.1} and \eqref{eq.KMP.3.16}, we conclude that uniformly for any $t \in \Lambda_T$ it holds
\beam \label{eq.KMP.3.17} \notag
&&J_1(x,\,m,\,t)= \sum_{n=1}^{m} \PP\left[ \sum_{i=1}^{n}  {\bf X}^{(i)}\,e^{-\xi(T_i)} \in x\,A \,,\;N(t)=n \right] \\[2mm]
&&\sim  \sum_{n=1}^{m}\sum_{i=1}^{n}  \PP\left[ {\bf X}^{(i)}\,e^{-\xi(T_i)} \in x\,A \,,\;N(t)=n \right]\\[2mm] \notag
&&\sim  \left(\sum_{n=1}^{\infty} - \sum_{n=m+1}^{\infty} \right)\sum_{i=1}^{n}  \PP\left[ {\bf X}^{(i)}\,e^{-\xi(T_i)} \in x\,A \,,\;N(t)=n \right] =:J_{11}(x,\,m,\,t)+J_{12}(x,\,m,\,t)\,.
\eeam

For the first term we find
\beam \label{eq.KMP.3.18} 
J_{11}(x,\,m,\,t)= \sum_{i=1}^{n}  \PP\left[ {\bf X}^{(i)}\,e^{-\xi(T_i)} \in x\,A \,,\;N(t)\geq i \right] =\int_0^t \PP[{\bf X}^{(i)}\,e^{-\xi(s)} \in x\,A ]\,\lambda(ds)\,.
\eeam

The second term, through \eqref{eq.KMP.3.10}, gives
\beam \label{eq.KMP.3.19}
&&J_{12}(x,\,m,\,t) \lesssim k^* \sum_{n=m+1}^{\infty} \sum_{i=1}^{n}  \PP\left[ {\bf X}^{(i)}\,e^{-\xi(T_i)} \in x\,A \right]\,\PP[N(t)=n] \\[2mm] \notag
&&=k^* \sum_{n=m+1}^{\infty} n  \PP\left[ {\bf X} \in x\,A \right]\,\PP[N(t)=n]=k^* \PP[{\bf X} \in x\,A ]\,\E\left[N(t)\,{\bf 1}_{\{N(t)> m\}}\right]\,,
\eeam
where the $k^*:=\vee_{i\in \bbn} k_2(i)$. Hence, by \eqref{eq.KMP.3.4} and \eqref{eq.KMP.3.19} we find
\beam \label{eq.KMP.3.20} \notag
&&\lim_{m\to \infty} \limsup_{\xto} \sup_{t \in \Lambda_T} \dfrac{J_{12}(x,\,m,\,t)}{\int_0^t \PP[{\bf X}^{(i)}\,e^{-\xi(s)} \in x\,A ]\,\lambda(ds)}  \\[2mm]
&&\leq \lim_{m\to \infty} \limsup_{\xto} \sup_{t \in \Lambda_T} \dfrac{k^*\,\PP[{\bf X} \in x\,A ]\,\E\left[N(t)\,{\bf 1}_{\{N(t)> m\}}\right]}{K_2\,\PP\left[ {\bf X} \in x\,A \right]\,\,\lambda(t)}=0\,,  
\eeam
where in last step we used again \cite[Lem. 3.2]{tang:2007}. From relations \eqref{eq.KMP.3.17}, \eqref{eq.KMP.3.18} and \eqref{eq.KMP.3.20} we conclude that 
\beam \label{eq.KMP.3.21} 
J_1(x,\,m,\,t) \sim \int_0^t \PP[{\bf X}^{(i)}\,e^{-\xi(s)} \in x\,A ]\,\lambda(ds)\,. 
\eeam
From \eqref{eq.KMP.3.1}, in combination with \eqref{eq.KMP.3.5} and \eqref{eq.KMP.3.21}, we get relation \eqref{eq.KMP.2.18} uniformly for any $t \in \Lambda_T$.
~\halmos

{\bf Proof of Corollary \ref{cor.KMP.2.1}.}~
At first, we notice that by Assumption \ref{ass.KMP.2.1} and the condition of bounded premium densities, we can obtain that
\beam \label{eq.KMP.3.22} 
0 \leq  \int_0^t \,e^{-\xi(y)} c_i(y)\,dy \leq M_i\, e^{C_1}\,T < \infty\,, 
\eeam
for any $i=1,\,\ldots,\,d$ and $t \in \Lambda_T$. From relation \eqref{eq.KMP.3.22}
and \cite[Lem. 4.3(d)]{samorodnitsky:sun:2016}, since $A=({\bf l} - L) \in \mathscr{R}$ we conclude that there exists some $u>0$, such that it holds
\beam \label{eq.KMP.3.23} 
(x+u)\,A \subsetneq x\,A \subsetneq (x-u)\,A \,. 
\eeam
From \eqref{eq.KMP.3.23} and Theorem \ref{th.KMP.2.1} in combination with class $\mathcal{L}_A$ properties, we obtain
\beam \label{eq.KMP.3.24} \notag 
&&\psi_{{\bf l},L}(x;\,t)\leq \PP\left[{\bf D}(t) \in x\,A + \int_0^s \,e^{-\xi(y)} c_i(y)\,dy\,, \;\exists \; s \in [0,\,t] \right] \leq \PP\left[{\bf D}(t) \in (x-u)\,A \right]\\[2mm] 
&& \sim \int_0^t \PP\left[{\bf X}\,e^{-\xi(s)} \in (x-u)\,A \right]\,\lambda(ds) \sim  \int_0^t \PP\left[{\bf X}\,e^{-\xi(s)} \in x\,A \right]\,\lambda(ds) \,,
\eeam
uniformly for any $t \in \Lambda_T$. From the other hand side, similar arguments provide the next relation
\beam \label{eq.KMP.3.25} \notag 
&&\psi_{{\bf l},L}(x;\,t)\geq \PP\left[{\bf D}(t)- \int_0^t \,e^{-\xi(y)} c_i(y)\,dy \in x\,A \right] \geq \PP\left[{\bf D}(t) \in (x+u)\,A \right]\\[2mm] 
&& \sim \int_0^t \PP\left[{\bf X}\,e^{-\xi(s)} \in (x+u)\,A \right]\,\lambda(ds) \sim  \int_0^t \PP\left[{\bf X}\,e^{-\xi(s)} \in x\,A \right]\,\lambda(ds) \,,
\eeam
uniformly for any $t \in \Lambda_T$. Therefore, by  \eqref{eq.KMP.3.24} and \eqref{eq.KMP.3.25}, we obtain  \eqref{eq.KMP.2.24} uniformly for any $t \in \Lambda_T$.
~\halmos

\subsection{Infinite time horizon}

Before the proof of Theorem \ref{th.KMP.4.1}, we establish a version of Breiman's lemma on set from family $\mathscr{R}$. We see that the next lemma serves to help the reader and not to generalized the Breiman's lemma in multivariate set up. Indeed, this way we can avoid the vague convergence and to apply directly the asymptotic expressions. For proper, independent and dependent multivariate versions of the Breiman's lemma we refer to \cite{basrak:davis:mikosch:2002} and \cite{fougeres:mercadier:2012} respectively.

\ble \label{lem.KMP.5.1}
Let $A \in \mathscr{R}$ be some fixed set. If ${\bf X}$, $\Theta$ are independent, with ${\bf X}$ to follow distribution $F \in MRV(\alpha,\,V,\,\mu)$, with $\alpha \in (0,\,\infty)$, and $\E[\Theta^{p} ]<\infty$, for some $p> \alpha$, then it holds
\beam \label{eq.KMP.3.26}
\PP[\Theta\,{\bf X} \in x\,A] \sim \E[\Theta^{\alpha} ]\,\PP[{\bf X} \in x\,A]\,,
\eeam
and furthermore, the distribution of $\Theta\,Y_A$ belongs to $\mathcal{R}_{-\alpha}$.
\ele

\pr~
At first, from the assumption of the finite moment of $\Theta$, applying the Breiman's lemma, see in \cite[Cor. 3.6]{cline:samorodnitsky:1994} or in \cite[Prop. 5.2(iv)]{leipus:siaulys:konstantinides:2023}, we obtain
\beam \label{eq.KMP.3.27}
\PP[\Theta\,{\bf X} \in x\,A] =\PP[\Theta\,Y_A > x] \sim \E[\Theta^{\alpha} ]\,\PP[Y_A > x] =\E[\Theta^{\alpha} ]\,\PP[{\bf X} \in x\,A] \,,
\eeam
which implies \eqref{eq.KMP.3.26}. Further, since $F_A \in \mathcal{R}_{-\alpha}$, through the closure property of the regular variation with respect to strong tail equivalence, see in \cite[Prop. 3.3(i)]{leipus:siaulys:konstantinides:2023}, we find that the distribution of $\Theta\,Y_A$ belongs to $\mathcal{R}_{-\alpha}$.
~\halmos

{\bf Proof of Theorem \ref{th.KMP.4.1}.}~
At first, from $F \in MRV(\alpha,\,V,\,\mu)$, follows, through Lemma \ref{lem.KMP.5.1}, that for any $\vep \in (0,\,1) $ there exists some large enough $x_0>0$, such that it holds
\beam \label{eq.KMP.3.28}
1-\vep \leq \dfrac {\PP[\Theta\,{\bf X} \in x\,A]}{\E[\Theta^{\alpha} ]\,\PP[{\bf X} \in x\,A] } \leq 1+\vep \,,
\eeam
for any $x\geq x_0$. Putting instead of $\Theta$ the $e^{-\xi(s)}$, because of Assumption \ref{ass.KMP.4.1}, we obtain for any $x\geq x_0$
\beam \label{eq.KMP.3.29}
\dfrac {\int_t^{\infty} \PP[{\bf X}\,e^{-\xi(s)} \in x\,A]\,\lambda(ds)}{\int_0^{\infty} \PP[{\bf X}\,e^{-\xi(s)} \in x\,A]\,\lambda(ds)} \leq \dfrac{(1+\vep)\PP[{\bf X} \in x\,A]\int_t^{\infty} \E\left[e^{-\alpha \,\xi(s)} \right]\,\lambda(ds)}{(1-\vep)\PP[{\bf X} \in x\,A]\int_0^{\infty} \E\left[e^{-\alpha \,\xi(s)} \right]\,\lambda(ds)} \,.
\eeam
From \eqref{eq.KMP.3.29}, for any $\delta>0$, we can find some $T_0 \in \Lambda$, such that the inequality
\beam \label{eq.KMP.3.30}
\int_{T_0}^{\infty} \PP[{\bf X}\,e^{-\xi(s)} \in x\,A]\,\lambda(ds)
\leq \delta\,\int_0^{T_0} \PP[{\bf X}\,e^{-\xi(s)} \in x\,A]\,\lambda(ds)\,,
\eeam
is true. By Theorem \ref{th.KMP.2.1} and relation \eqref{eq.KMP.3.30}, we conclude that it holds
\beam \label{eq.KMP.3.31} \notag
&&\PP[{\bf D}(t) \in x\,A]\geq \PP[{\bf D}(T_0) \in x\,A]\sim \int_0^{T_0} \PP[{\bf X}\,e^{-\xi(s)} \in x\,A]\,\lambda(ds)\\[2mm] \notag
&&\geq \left(\int_0^{t} -\int_{T_0}^{\infty}\right)\PP[{\bf X}\,e^{-\xi(s)} \in x\,A]\,\lambda(ds) \geq (1-\delta)\int_0^{t}\PP[{\bf X}\,e^{-\xi(s)} \in x\,A]\,\lambda(ds) \\[2mm]
&&\geq (1-\delta)\,(1 - \vep)\PP[{\bf X} \in x\,A]\int_0^{t}\E\left[e^{-\alpha \,\xi(s)} \right]\,\lambda(ds) \,,
\eeam
uniformly for any $t \in [T_0,\,\infty)$, where in last step we used relation \eqref{eq.KMP.3.28}.

From the other side, from \cite[Prop. 2.4]{konstantinides:passalidis:2025g} in second step and from \cite[Th. 3.3]{chen:yuen:2009} in the third step, due to Assumption \ref{ass.KMP.4.1} and Assumption \ref{ass.KMP.2.2}, relations \eqref{eq.KMP.3.30} and \eqref{eq.KMP.3.28} we obtain
\beam \label{eq.KMP.3.32} \notag
&&\PP[{\bf D}(t) \in x\,A]\leq \PP[{\bf D}(\infty) \in x\,A] \leq  \PP\left[\sum_{i=1}^{\infty} Y_A^{(i)}\,e^{-\xi(T_i)} > x\right]\sim \sum_{i=1}^{\infty}\PP\left[Y_A^{(i)}\,e^{-\xi(T_i)} > x\right]\\[2mm]
&&=\int_0^{\infty}\PP[{\bf X}\,e^{-\xi(s)} \in x\,A]\,\lambda(ds)\leq \left(\int_0^{t} +\int_{T_0}^{\infty}\right)\PP[{\bf X}\,e^{-\xi(s)} \in x\,A]\,\lambda(ds)  \\[2mm] \notag
&&\leq (1+\delta)\,(1 + \vep)\PP[{\bf X} \in x\,A]\int_0^{t}\E\left[e^{-\alpha \,\xi(s)} \right]\,\lambda(ds) \,,
\eeam
uniformly for any $t \in [T_0,\,\infty)$. Therefore, by \eqref{eq.KMP.3.31} and \eqref{eq.KMP.3.32}, choosing $M^-:=(1-\vep)\,(1-\delta)$ and $M^+:=(1+\vep)\,(1+\delta)$, given the arbitrarily choice of $\vep$  and $\delta$, we get by tending $\vep\downarrow 0$ and $\delta \downarrow 0$ the convergences  $M^-\uparrow 1$ and $M^+\downarrow 1 $ to reach the desired result for any $t \in \Lambda \cap (T_0,\,\infty]$.

By Theorem \ref{th.KMP.2.1} and via Lemma \ref{lem.KMP.5.1}, relation \eqref{eq.KMP.3.28} we find \eqref{eq.KMP.4.3} to hold uniformly for any $t \in \Lambda$, we omit the analytic presentation, as it follows similar steps. 
~\halmos

{\bf Proof of Corollary \ref{cor.KMP.4.1}.}~
For the lower bound of \eqref{eq.KMP.4.5}, we use the assumption that the integrals of total premiums of the $d$ lines of business are finite and through the inclusion of the set of \eqref{eq.KMP.3.23} we find for $u>0$ it holds
\beam \label{eq.KMP.3.33}
&&\psi_{{\bf l},L}(x,\,t)=\PP\left[{\bf D}(t)- \int_0^s \,e^{-\xi(y)} {\bf c}(y)\,dy \in x\,A\,,\;\exists \;s\in (0,\,t] \right] \\[2mm]  \notag
&& \geq \PP\left[{\bf D}(t)- \int_0^t \,e^{-\xi(y)} {\bf c}(y)\,dy \in x\,A  \right]\geq \PP\left[{\bf D}(t) \in (x+u)\,A \right]\\[2mm] \notag
&&  \sim \PP[{\bf X} \in (x+u)\,A]\,\int_0^t \E\left[e^{-\alpha \,\xi(s)} \right]\,\lambda(ds) \sim   \PP[{\bf X} \in x\,A]\,\int_0^t \E\left[e^{-\alpha \,\xi(s)} \right]\,\lambda(ds)  \,,
\eeam
uniformly for any $t \in \Lambda$, where in the last step we used the inclusion property $MRV(\alpha,\,V,\,\mu) \subsetneq \mathcal{L}_A$, see Section 2. From the other side, through relation \eqref{eq.KMP.3.23} and by similar arguments with \eqref{eq.KMP.3.32}, we get that it holds
\beam \label{eq.KMP.3.34} \notag
&&\psi_{{\bf l},L}(x,\,t)\leq \PP\left[{\bf D}(\infty)- \int_0^s \,e^{-\xi(y)} {\bf c}(y)\,dy \in x\,A\,,\;\exists \;s >0 \right] 
\leq \PP\left[{\bf D}(\infty)\in (x-u)\,A  \right]\\[2mm] 
&& \leq (1+\delta)\,(1+\vep)\, \PP[{\bf X} \in (x-u)\,A]\,\int_0^t \E\left[e^{-\alpha \,\xi(s)} \right]\,\lambda(ds) \\[2mm] \notag
&& \sim (1+\delta)\,(1+\vep)\, \PP[{\bf X} \in x\,A]\,\int_0^t \E\left[e^{-\alpha \,\xi(s)} \right]\,\lambda(ds)  \,,
\eeam
uniformly for any $t \in [T_0,\,\infty]$ and additionally, by relation \eqref{eq.KMP.3.24}, we reach the desired upper bound, for any $t \in \Lambda \cap [0,\,T_0]$. Hence, taking into account \eqref{eq.KMP.3.34}  and \eqref{eq.KMP.3.33},  we conclude that \eqref{eq.KMP.4.5} holds uniformly for any $t \in \Lambda$.
~\halmos

\noindent \textbf{Acknowledgments.} 
This work is dedicated to the memory of V.V. Kalashnikov (1942 - 2001) and R. Norberg (1945 - 2017).

\end{document}